# Optimal Design of Distributed Energy Systems Considering the Impacts on Electrical Power Networks

*Ishanki De Mel[a], Oleksiy V. Klymenko[a], Michael Short[a*]*

[a]Department of Chemical and Process Engineering, University of Surrey, Guildford, United Kingdom GU2 7XH

*m.short@surrey.ac.uk

## Abstract

Distributed energy systems (DES) have the potential to minimise costly network upgrades while increasing the proportion of renewable energy generation in the electrical grid, when properly designed. In contrast, poorly designed DES can accelerate the degeneration of existing network infrastructure. Most optimisation models used to design grid-connected DES have oversimplified or excluded constraints associated with alternating current (AC) power flow, as the latter has been studied in a standalone class of models known as Optimal Power Flow (OPF). A small subset of models, labelled DES-OPF models, have attempted to combine these two independent frameworks. However, the impacts of using a DES-OPF framework on the resulting designs, as opposed to a conventional DES framework without AC power flow, have not been studied in previous work. This study aims to shed light on these impacts by proposing a bi-level method to solve the computationally-expensive DES-OPF framework, and simultaneously comparing results to a baseline MILP model that utilises direct current (DC) approximations in place of AC OPF, as found abundantly in literature. Two test cases of varying scale are employed to test the frameworks and compare resulting designs. The practical feasibility of the designs is assessed, based on whether the designs can mitigate network violations and energy wastage during standard operation. Results demonstrate that the baseline MILP underestimates total costs due to its inability to detect current and voltage violations, resulting in a 37% increase for one case study when tested with the DES-OPF framework. Major implications on battery capacity are also observed, primarily to manage higher levels of renewable energy curtailment, which emphasise the need to use DES-OPF frameworks when designing grid-connected DES.

**Keywords:** Distributed Energy, Optimal Power Flow, Nonlinear Programming, Mixed-integer programming, Optimization.

## 1. Introduction

Ageing power networks around the world are being transformed by the increasing integration of renewable and distributed energy resources (DERs). Previously unidirectional distribution networks have been converted to bidirectional active networks through collections of small-scale energy generation and storage technologies situated in proximity to the consumer. These are known as Distributed Energy Systems (DES) and consist of mostly renewable and low-carbon technologies. The term 'DES' has become synonymous with microgrids, energy hubs, distributed generation (DG), and smart local energy systems (De Mel et al., 2022). In most European countries, such as the United Kingdom, the technologies within a DES are often owned by the consumers themselves, who benefit from cost-reducing tariffs attached to low-carbon energy generation (Ofgem, 2020; Ogfem, 2018). When designed appropriately with operational and network-related constraints in mind, DES have the potential to lower energy losses and carbon emissions related to power generation, and reduce the need for costly network reinforcements (Prakash and Khatod, 2016; UK Power Networks, 2018). However, if not designed and operated with care, these systems can potentially amplify voltage unbalance, increase power losses (such as via renewable energy curtailment), and subsequently lead to large economic losses (Kuang et al., 2011). Due to these implications, the design of DES is often modelled as a mathematical optimisation problem. This can ensure the suitability of selected technologies over the system lifetime, with respect to relevant constraints and objective functions, such as the minimisation of annualised costs and/or carbon emissions (Mavromatidis et al., 2019).

Optimisation models account for operational aspects such as average demand and weather profiles, and to aid detailed analyses of how the system responds to dynamic conditions and requirements.

DES design is most commonly represented using mixed-integer formulations, where discrete decisions are used to choose the location of installation, the type of technology to be installed, and different states of operation. Continuous variables are used to describe energy flows, voltages, costs, etc. Although nonlinear formulations often more realistically depict the behaviour of certain technologies and underlying networks, most DES models have been formulated and solved as Mixed-Integer Linear Programming (MILP) problems (Theo et al., 2017). Despite most DES operating in grid-connected mode, these MILP models oversimplify or exclude nonlinearities associated with the underlying electrical network. For example, recent studies that investigate the use of advanced features within DES design frameworks, such as bidirectional power flow in peer-to-peer trading (Meena et al., 2019) with blockchain protocols (Li et al., 2019), utilising batteries to minimise stresses placed on the grid (Liu et al., 2020), and the integration of hydrogen storage (Xiang et al., 2021), have all overlooked the inclusion of constraints related to the distribution network. Distribution networks are part of the larger electrical grid, which transmits and distributes alternating current (AC) power. They are utilised by DES for active and reactive power flows, while the network also acts as both source and sink when demand within the DES does not match supply. Most DES design models only consider an active power balance, and exclude calculations of reactive power, voltages, and currents. These studies often employ the DC approximation (Sass et al., 2020), which is widely considered as a simplified alternative to AC power flow. However, this does not accurately represent AC power flow as it calculates only the active power flow at each branch of the network, assuming that voltage magnitudes remain at nominal values and neglecting differences between voltage angles (Baker, 2020; Purchala et al., 2005).

In a recent literature review (De Mel et al., 2022), we identified a small subset of models labelled 'DES-OPF' that incorporated detailed representations of AC power flow, traditionally found in a separate class of models known as Optimal Power Flow (OPF), to aid better decision-making with respect to grid-connected DES. OPF models have been used to study the allocation of distributed generators in transmission networks and contain detailed nonlinear and nonconvex constraints related to active and reactive power flows, voltages, and currents. Unlike DES models, these do not include constraints describing the operation of distributed generators and storage technologies. These models primarily use balanced power flow equations, as the powers, voltages, and currents are assumed to be balanced across all phases. Despite the additional complexity introduced by AC OPF, DES and OPF can no longer be treated as standalone problems, as DERs installed within a DES have the potential to disrupt the normal functioning of the underlying network. Including OPF in the earlier stage of DES planning and implementation, such as in design models, has the potential to eliminate issues that may arise at the operational level, such as elevated network unbalances and reduced power quality. This would result in more robust DES that can operate safely within existing distribution networks. Therefore, the need to achieve an accuracy-complexity balance when formulating and solving DES-OPF optimisation models is an active area of research.

The literature review (De Mel et al., 2022) further highlighted that few studies exist on DES design, as most DES-OPF frameworks focus on operation or scheduling decisions based on the assumption that the design has been predetermined. Morvaj et al. (2016) and Mashayekh et al. (2017) have proposed linearised AC power flow equations that can be implemented within MILP design formulations. A post-optimisation check using MATPOWER (Zimmerman et al., 2011), a popular power flow simulation tool for balanced AC networks, and a Genetic Algorithm (GA) design model have also been employed by Morvaj et al. (2016). The former is a key weakness of the framework, as the AC nonconvex equations have no influence on the DES design, thus merely serving as a numerical check. Linearised power flow equations have been proposed in both these studies to remedy this and are shown to have high accuracy when compared to the results original nonconvex formulations. However, the implications of these linearisations on the scalability of this framework have not been studied. Morvaj et al. (2016)

also report many current violations within the distribution network as a result of DES implementation, highlighting the need for further research in this area. In contrast, Sfikas et al. (2015) and Rezaee Jordehi et al. (2021) have utilised nonlinear AC power flow equations as opposed to linear approximations. The former study analyses the suitability of the design in both grid-connected and islanded modes, but have excluded the use of binary decisions to solve the model using Sequential Quadratic Programming (SQP). This results in a more simplistic representation of the technologies and their operation. In comparison, the latter study (Rezaee Jordehi et al., 2021) utilises a more detailed Mixed-Integer Nonlinear Programming (MINLP) framework, where results demonstrate how the placement of a battery swapping station can impact power losses and voltages in the network. MINLP formulations typically tend to be harder to solve and more computationally expensive, often resulting in intractable models when scaled up.

All of the studies mentioned above use balanced AC power flow equations where only a single phase is analysed, as done in OPF models for transmission networks. In reality, small-scale technologies within DES are connected to low-voltage distribution networks which tend to be highly unbalanced due to its radial structure and varying cumulative loads connected to each phase (Ma et al., 2020). Dunham et al. (2020) attempt to remedy this by incorporating an iterative fixed-point linear approximation of multiphase AC power flow (Bernstein and Dall'anese, 2017) in the DES design model. Balanced OPF offers a significant advantage over multiphase OPF due to the reduced complexity, resulting in much more tractable models, especially when integer decisions are also involved for DES placement and operation. Furthermore, the linearisation of multiphase power flow requires feasible initial conditions or fixed points, which may need to be determined or iteratively verified by solving the full power flow formulations using power flow tools (Bernstein et al., 2018) such as OpenDSS (Electric Power Research Institute, 2020).

## 1.1. Contributions of this study

While all of the DES-OPF studies discussed above indicate the need for a combined framework that incorporates AC power flow, there has been a significant lack of comparison to DES MILP frameworks that use DC approximations (henceforth labelled *DES baseline models*), which are abundant in literature. There is little knowledge or quantitative analysis on how DES designs are impacted when OPF formulations are included, and how operational strategies may also be impacted as a result. Furthermore, existing DES-OPF frameworks either employ linear approximations with an increased number of variables when compared to the original formulations, or use MINLP frameworks that may not be suitable for test cases involving greater number of nodes and buses. The use of metaheuristic techniques and post-optimisation checks cannot guarantee local nor global optimality, which highlights the need for a mathematical programming approach that can produce locally optimal and feasible results that can also be implemented in practice. We argue that the feasibility of the solutions is more important than global optimality in this instance, as the solutions of these models should guarantee the synergistic operation of the DES within the existing distribution network.

To address these gaps, we propose a bi-level framework where OPF constraints can influence DES design decisions without the use of post-optimisation checks, linear approximations, or metaheuristics. We then test the following hypothesis:

*Designs from DES baseline models may not be practically feasible when tested with a DES-OPF framework.*

As most optimisation models for DES design aim to introduce or increase the proportion of renewable energy utilised by the consumers while minimising costs, the term 'practically feasible' in the hypothesis refers to the assessment of the design post-optimisation based on the following criteria:

1. Can the design respect network constraints during baseline operation, and thus be implemented?

2. Does the design prevent or minimise energy curtailment or wastage during baseline operation?

The criteria focus on baseline operation as this represents the operation of DES in the absence of serious disturbances or anomalies. It may not be possible to avoid network violations or energy curtailment in the presence of disturbances.

To test the hypothesis, we use two test cases, including a modified version of the IEEE EU LV network, to generate results from the proposed bi-level method and compare to those from a DES baseline model. Results are also compared with an equivalent MINLP framework to assess computational expense and tractability.

With the increasingly multi-faceted nature of DES optimisation problems, the potential implications of a DES design not being compatible with existing infrastructure cannot be overlooked. Such designs can decrease efficiency while incurring additional costs to both consumers and network operators, straying from the true objectives of implementing a DES to reduce environmental impacts and minimise costs for all parties.

## 2. Methodology

Figure 1 summarises the frameworks employed in this study, where i) MILP refers to the DES baseline model that uses the DC approximation, ii) Bi-level (BL) refers to the bi-level method incorporating nonlinear AC power flow proposed in this study, and iii) MINLP is the mixed-integer nonlinear programming model which is solved simultaneously using commercial solvers.

Focusing on the bi-level method, referred to as BL henceforth, Figure 1 highlights that the entire MILP framework representing the DES baseline model is solved first, under Level 1, to obtain the binary topology for a globally optimal DES design excluding detailed network constraints. The binary topology contains both design and operational binary variables, and these are fixed before the nonlinear OPF constraints are introduced in Level 2. This allows the combined DES-OPF framework to be solved using nonlinear programming (NLP), as only continuous variables (including key design variables such as capacities of units) are retained in the formulation.

The proposed bi-level method offers some notable advantages compared to previous DES-OPF frameworks. Firstly, the OPF constraints have the ability to influence design variables and revise the initial design proposed by the baseline MILP, such that the network constraints are not violated during normal operation. This is superior to performing post-optimisation numerical checks, as it eliminates the need to solve optimisation models iteratively, in conjunction with external power flow tools, to achieve convergence. Secondly, solving the MILP containing the DES constraints first and fixing the binary topology allows the overall DES-OPF framework to be solved as a large-scale NLP model, which enables the use of reputable NLP solvers such as CONOPT (Drud, 1985) and IPOPT (Wächter and Biegler, 2006) to obtain faster results and guaranteed locally optimal solutions. This is preferred over solving the formulation as a costly MINLP or via less reliable metaheuristic optimisation approaches that may require many tuning parameters. Thus, the framework strives to achieve an accuracy-complexity balance, which is essential for solving larger DES. Finally, as portrayed in Figure 1, the overall framework facilitates the comparison with alternative solving methods, such as the uncoupling of the nonlinear OPF formulation to solve the model as a baseline MILP, or the retention of binary variables to solve the combined DES-OPF formulation as an MINLP.

To test whether the baseline MILP design can achieve a feasible operational schedule with respect to the network as well, a modified version of BL (labelled *MILP Check*) is used, where design variables are also fixed alongside the binary topology. This restricts the degrees of freedom in the DES-OPF model to continuous operational variables only, while the design proposed by the DES baseline must remain unchanged.

In line with the bi-level framework, the complete DES Design – OPF formulation can be broken down into two subsections, Level 1: DES Design (Section **Error! Reference source not found.**) and Level 2: OPF (Section **Error! Reference source not found.**), which are subsequently linked by a number of linking constraints (Section 2.3) in the combined framework. Note that the objective function, applicable to both DES and DES-OPF models, is presented in Section **Error! Reference source not found.**.

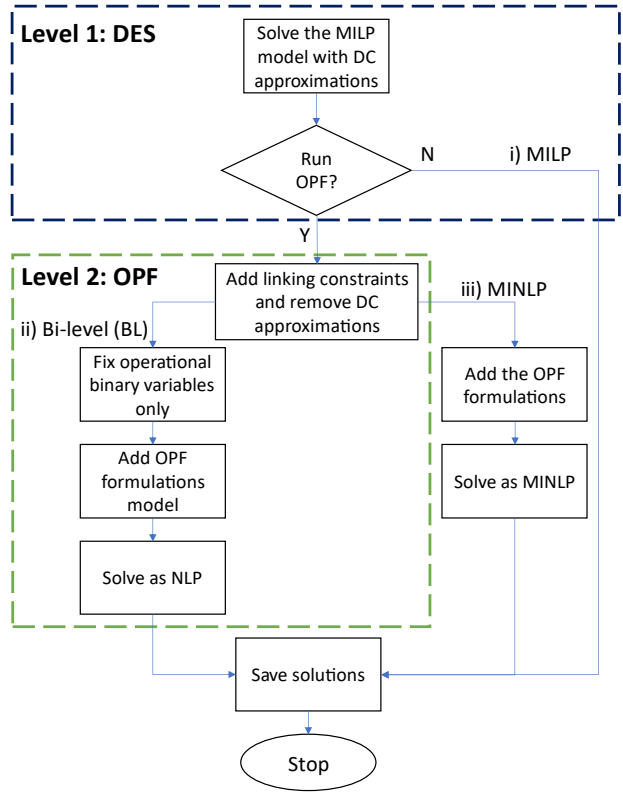

Figure 1. The DES-OPF modelling approaches employed in this study, where i) an MILP can be solved with DC approximations for comparison, ii) BL is the bi-level method proposed in this study, which allows the OPF constraints to influence both design and operational variables while keeping binary variables fixed, iii) an MINLP model where OPF constraints can influence all design and operational aspects of the DES model.

## 2.1. Level 1: DES Design

An MILP formulation is utilised as the DES baseline model, based on the works of Ren and Gao (2010), Mehleri et al. (2012), and Mariaud et al. (2017). Several modifications have been made to add flexibility and enable facile implementation when using the proposed DES – OPF methods, which are described in this section.

The base model has been formulated such that it can be solved independently as a seasonal model, if required, or across any number of seasons $s \in S$. A common objective function to minimise the total annualised cost $TAC$, which is applicable to both DES and DES-OPF frameworks, is used across all seasons:

$$\min \ TAC = \sum_{G \in DER} C_G^{INV} + \sum_{s \in S} \left( C_s^{grid} + \left( \sum_{G \in DER} C_{s,G}^{OM} \right) - CI_s \right) \tag{1}$$

The total annualised cost is the sum of investment costs $C_G^{INV}$ for the generalised set of DERs $G \in DER$, seasonal costs for purchasing electricity from the grid $C_s^{grid}$ and operation of DERs $C_{s,G}^{OM}$, and seasonal income $CI_s$. Linking constraints across seasons are also utilised to ensure that the capacity of each DER $Cap_{i,s}^G$ installed at each residence or building $i \in \mathbf{I}$ does not vary from season to season:

$$If \ s > 1 \quad Cap_{i,s}^{G} = Cap_{i,s-1}^{G} \quad \forall i \in \mathbf{I} \tag{2}$$

Thus, the annualised investment costs for each DER can be calculated with respect to the selected capacity, capacity cost $\mathrm{CC}_{G}^{\mathrm{INV}}$(£/kW), and capital recovery factor CRF which takes the total lifetime of the DES into account:

$$C_{G}^{INV} = Cap_{i,s}^{G} * \mathrm{CC}_{G}^{\mathrm{INV}} * \mathrm{CRF} \quad \forall G \in \mathrm{DER} \tag{3}$$

The constraints included in each seasonal model are discussed below. Note that the subscript $s \in S$, representing seasons, has now been removed from the formulation. Electrical power flows are denoted by $E$, while thermal power flows are denoted by $H$, with both calculated in kilowatts (kW).

Electricity demand $\mathrm{E}_{i,t}^{\mathrm{load}}$ at each building $i \in \mathbf{I}$ and timepoint $t \in \mathbf{T}$ can be satisfied using the equation below:

$$\mathrm{E}_{i,t}^{\mathrm{load}} = E_{i,t}^{grid} + E_{i,t}^{G} + \sum_{c} E_{i,t,c}^{disch} \quad \forall i \in \mathbf{I}, t \in \mathbf{T} \tag{4}$$

Where $E_{i,t}^{grid}$ is the purchased power, $E_{i,t}^{G}$ is the power generated by the electricity-generating distributed generators $G \in DER$, and $E_{i,t,c}^{disch}$ is the power discharged by the batteries $c \in \mathbf{C}$ to satisfy consumer demand.

The power generated onsite must be used first to satisfy demand, allowing customers to sell excess power only, $E_{i,t}^{sold}$. This results in a bilinear constraint, reformulated using binary $X_{i,t}$ and a big-M constraint as shown below (Mehleri et al., 2012):

$$E_{i,t}^{grid} \le \mathrm{E}_{i,t}^{\mathrm{load}} * \left(1 - X_{i,t}\right) \quad \forall i \in \mathbf{I}, t \in \mathbf{T} \tag{5}$$

$$E_{i,t}^{sold} \le \mathrm{M} * X_{i,t} \quad \forall i \in \mathbf{I}, t \in \mathbf{T} \tag{6}$$

Heat demand $\mathrm{H}_{i,t}^{\mathrm{load}}$ is met using the thermal power generated by heat-generating DERs, $H_{i,t}^{G}$:

$$\mathrm{H}_{i,t}^{\mathrm{load}} = H_{i,t}^{G} \quad \forall i \in \mathbf{I}, t \in \mathbf{T} \tag{7}$$

Three of the most prevalent DERs implemented in residential settings have been included in this formulation, which are solar photovoltaics (PVs), batteries, and boilers. Note that additional technologies such as heating/cooling networks and other types of DERs have not been included in the formulation as the aim of this work is to study how electricity-generating technologies can impact the electricity distribution network. The technology-specific superscripts are used instead of superscript $G$ henceforth.

For solar PVs specifically, solar irradiance available at each timepoint $\mathrm{Irr}_{t}$ (kW/m²), panel efficiency $\eta^{\mathrm{panel}}$, number of panels installed at each building $N_{i}^{Panel}$, and area of each panel $\mathrm{A}^{\mathrm{Panel}}$ (m²) limits the amount of total power produced:

$$E_{i,t}^{PV,used} + E_{i,t}^{PV,sold} + \sum_{c} E_{i,t,c}^{PV,charge} \le N_{i}^{Panel} * \mathrm{A}^{\mathrm{Panel}} * \mathrm{Irr}_{t} * \eta^{\mathrm{panel}} \quad \forall i \in \mathbf{I}, t \in \mathbf{T} \tag{8}$$

The total solar power generated is broken down to represent the power consumed by the building to meet demand $E_{i,t}^{PV,used}$, the excess power sold $E_{i,t}^{PV,sold}$, and power used to charge batteries $E_{i,t,c}^{PV,charge}$.

The maximum available roof area $\mathrm{A}_{i}^{\mathrm{roof}}$(m²) at each building limits the number of panels that can be installed:

$$N_{i}^{Panel} * \mathrm{A}^{\mathrm{Panel}} \le \mathrm{A}_{i}^{\mathrm{roof}} \quad \forall i \in \mathbf{I} \tag{9}$$

To calculate the total operational and maintenance costs of installing PVs $C^{OM,PV}$, the fixed yearly operational cost parameter $C^{OM,f,PV}$ (£/kW-year) is accounted for with respect to the total capacity of panels installed, where $Cp^{Panel}$ is the rate capacity of one solar panel:

$$C^{OM,PV} = \sum_i N_i^{Panel} * C^{OM,f,PV} * Cp^{Panel} * \frac{1}{365} * N^{days} \tag{10}$$

As boilers are utilised to generate heat $H_{i,t}^b$, the maximum capacity required to satisfy the demand is determined by the maximum heat generated $H_i^{b,max}$ (which is representative of the total capacity of the boiler $Cap_i^b$ required for Eqs. (2) and (3)):

$$H_i^{b,max} \geq H_{i,t}^b \quad \forall i \in \mathbf{I}, t \in \mathbf{T} \tag{11}$$

The total operational cost of boilers within the DES $C^{OM,b}$ is a function of the fuel price $C^{gas}$ (£/kWh) and thermal efficiency $\eta^b$:

$$C^{OM,b} = \sum_{i \in I, t \in T} H_{i,t}^b * \Delta t * C^{gas} / \eta^b \tag{12}$$

The capacity of the Lithium-ion batteries installed at each building $Cap_{i,c}{}^{batt}$ (kWh) can be determined by the volume installed $V_{i,c}$ (m³) and the volumetric energy density of the battery $VED_c$ (kWh/m³) (Mariaud et al., 2017):

$$Cap_{i,c}{}^{batt} = V_{i,c} * VED_c \quad \forall i \in \mathbf{I}, c \in \mathbf{C} \tag{13}$$

The volume of the battery installed cannot exceed the maximum available volume at each building $VA_i$:

$$\sum_c V_{i,c} \leq VA_i \quad \forall i \in \mathbf{I} \tag{14}$$

If multiple battery-options are presented (such as Li-ion, Sodium-Sulphur, etc.), a binary variable $W_{i,t}$ can be used within a big-M constraint to determine which type is installed (note that the options in this study are limited to Li-ion):

$$\sum_c W_{i,c} \leq 1 \quad \forall i \in \mathbf{I} \tag{15}$$

$$Cap_{i,c}{}^{batt} \leq M * W_{i,c} \quad \forall i \in \mathbf{I}, c \in \tag{16}$$

The maximum state of charge $SoC_c^{max}$ and depth of discharge $DoD_c^{max}$ allowed limits the amount of energy $\mathcal{E}_{i,t,c}^{stored}$ (kWh) a battery can store:

$$\mathcal{E}_{i,t,c}^{stored} \leq Cap_{i,c}{}^{batt} * SoC_c^{max} \quad \forall i \in \mathbf{I}, t \in \mathbf{T}, c \in \mathbf{C} \tag{17}$$

$$\mathcal{E}_{i,t,c}^{stored} \geq Cap_{i,c}{}^{batt} * (1 - DoD_c^{max}) \quad \forall i \in \mathbf{I}, t \in \mathbf{T}, c \in \mathbf{C} \tag{18}$$

The amount of energy stored in the battery is governed by the charging power $E_{i,t,c}^{ch}$, discharging power $E_{i,t,c}^{disch}$, and respective charging and discharging efficiencies, $\eta^{ch}$ and $\eta^{disch}$:

$$if\ t = start:\ \mathcal{E}_{i,t,c}^{stored} = \left(E_{i,t,c}^{ch} * \eta^{ch} * \Delta t\right) - \frac{E_{i,t,c}^{disch} * \Delta t}{\eta^{disch}} \quad \forall i \in \mathbf{I}, t \in \mathbf{T}, c \in \mathbf{C} \tag{19}$$

$$else:\ \mathcal{E}_{i,t,c}^{stored} = \mathcal{E}_{i,t-1,c}^{stored} + \left(E_{i,t,c}^{ch} * \eta^{ch} * \Delta t\right) - \frac{E_{i,t,c}^{disch} * \Delta t}{\eta^{disch}} \quad \forall i \in \mathbf{I}, t \in \mathbf{T}, c \in \mathbf{C} \tag{20}$$

Note that at the first timepoint is denoted as $t = start$. As no previous stored energy $\mathcal{E}_{i,t-1,c}^{stored}$ exists at this timepoint, an if-else statement is used, such that subsequent timepoints include the energy storage level at the previous timepoint.

To prevent the battery from discharging more energy than it already contains in the previous timepoint, a logical condition is presented:

$$If \ t > 1 \quad E_{i,t,c}^{disch} * \frac{\Delta t}{\eta^{disch}} \leq \mathcal{E}_{i,t-1,c}^{stored} \quad \forall i \in \mathbf{I}, t \in \mathbf{T}, c \in \mathbf{C} \tag{21}$$

The storage levels at the beginning and end of the temporal horizon, denoted by $t = start$ and $t = end$, are also fixed, to ensure that there are no unaccounted differences in the storage between each day of the respective season:

$$\mathcal{E}_{i,t=start,c}^{stored} = \mathcal{E}_{i,t=end,c}^{stored} \quad \forall i \in \mathbf{I}, c \in \mathbf{C} \tag{22}$$

The batteries can be charged using either PV power $E_{i,t,c}^{ch,PV}$ or power purchased from the grid $E_{i,t,c}^{ch,grid}$:

$$E_{i,t,c}^{ch} = E_{i,t,c}^{ch,grid} + E_{i,t,c}^{ch,PV} \quad \forall i \in \mathbf{I}, t \in \mathbf{T}, c \in \mathbf{C} \tag{23}$$

To prevent the battery from cycling unrealistically, scalar upper bounds $S_1$ and $S_2$ have been placed on charging and discharging at each timepoint:

$$E_{i,t,c}^{ch} * \Delta t * \eta^{ch} \leq Cap_{i,c}^{batt} * S_1 \quad \forall i \in \mathbf{I}, t \in \mathbf{T}, c \in \mathbf{C} \tag{24}$$

$$E_{i,t,c}^{disch} * \frac{\Delta t}{\eta^{disch}} \leq Cap_{i,c}^{batt} * S_2 \quad \forall i \in \mathbf{I}, t \in \mathbf{T}, c \in \mathbf{C} \tag{25}$$

Batteries are not allowed to charge and discharge at the same time, and therefore a binary variable $Q_{i,t,c}$ is used to indicate when the battery is charging:

$$E_{i,t,c}^{ch} \leq \mathrm{M} * Q_{i,t,c} \quad \forall i \in \mathbf{I}, t \in \mathbf{T}, c \in \mathbf{C} \tag{26}$$

$$E_{i,t,c}^{disch} \leq \mathrm{M} * (1 - Q_{i,t,c}) \quad \forall i \in \mathbf{I}, t \in \mathbf{T}, c \in \mathbf{C} \tag{27}$$

Note that big-M constraints are once again utilised to prevent bilinear terms from introducing nonlinearity into the linear DES model.

The total operational cost of batteries installed within the DES $C^{OM,Batt}$, is calculated using a fixed yearly operational cost $\mathrm{P}^{OM,batt}$(£/kW-year):

$$C^{OM,Batt} = \sum_i \frac{\sum_c Cap_{i,c}^{batt} * \mathrm{P}_c^{OM,batt}}{\Delta t} * \frac{1}{365} * N^{days} \tag{28}$$

The total purchasing cost of electricity for the whole DES $C^{grid}$ is calculated using the price of purchasing electricity $\mathrm{P}_t^{grid}$ (£/kWh), which can vary for day and night under Economy 7 tariffs (The Consumer Council, 2018):

$$C^{grid} = \sum_{i,t}(E_{i,t}^{grid} + \sum_c E_{i,t,c}^{ch,grid}) * \Delta t * \mathrm{P}_t^{grid} * \mathrm{N}^{days} \tag{29}$$

The total export income from selling excess power $CI^{export}$ can be calculated using the tariff prices offered by the Feed-In-Tariff (FIT) scheme (Ofgem, 2019) $\mathrm{P}^{tariff}$ (£/kWh), which is further described in Section 3:

$$CI^{export} = \sum_{i,t} E_{i,t}^{PV,sold} * \Delta t * \mathrm{P}^{tariff} * \mathrm{N}^{days} \tag{30}$$

The FIT also rewards consumers for generating renewable energy. The total generation income $CI^{gen}$ can be calculated using the generation-specific tariff $\mathrm{P}^{gen,tariff}$:

$$CI^{gen} = \sum_{i,t}(E_{i,t}^{PV,used} + E_{i,t}^{PV,sold} + \sum_c E_{i,t,c}^{charge}) * \Delta t * \mathrm{P}^{gen,tariff} * \mathrm{N}^{days} \tag{31}$$

The linear DC power flow approximation solved in Level 1 is presented below (Frank and Rebennack, 2016):

$$P_{n,t} \approx \sum_{m=1}^{N} \mathrm{B_{nm}}(\theta_{n,t} - \theta_{m,t}) \quad \forall n \in \mathbf{N}, t \in \mathbf{T}$$

Where $P_{n,t}$ is the active power at node $n$, $\mathrm{B_{nm}}$ is the susceptance at branch $(n, m)$ and $\theta_{n,t}$ is the voltage angle. The calculation for $\mathrm{B_{nm}}$ is explained further in Section 2.2 Eq. (12). The same linear linking constraint provided in Section 2.3 Eq. (22) for active power is utilised.

## 2.2. Level 2: OPF

The bus injection model presented in Frank and Rebennack (2016) is utilised for the OPF formulation, with an additional index representing time $t \in \mathbf{T}$. Complex voltage is denoted in polar form as $V \angle \theta^o$, where the magnitude of voltage is denoted as $V$ and the angle of the voltage is denoted as $\theta$. This formulation assumes that the low voltage network is balanced, and therefore considers only one phase in the formulation. The per unit system is also utilised in the OPF formulation, such that voltages, powers, and currents are made dimensionless or 'per unit', using balanced base values for voltage $\mathrm{V_{base}}$ and apparent power $\mathrm{S_{base}}$. This has several advantages, including ease of formulation and improvement of numerical stability (Frank and Rebennack, 2016). Figure 2 presents a line diagram of a small distribution network and introduces the set notations corresponding to the elements comprising the network. This notation is used in the OPF formulation below.

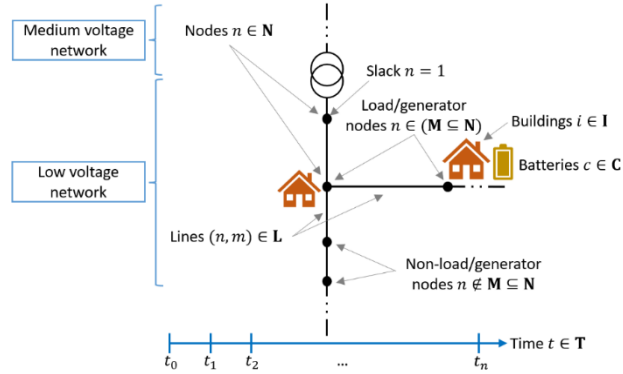

Figure 2. A graphical representation of the distribution network and the set notation used in the OPF formulation.

Considering a node $n \in \mathbf{N}$ on branch or line $(n, m) \in \mathbf{L}$, the active power $P_{n,t}$ and reactive power $Q_{n,t}$ balances can be described as a function of polar voltage, as shown in the nonconvex equality constraints below:

$$P_{n,t} = V_{n,t} \sum_{m=1}^{N} V_{m,t}((\mathrm{G_{nm}}\cos(\theta_{n,t} - \theta_{m,t}) + \mathrm{B_{nm}}\sin(\theta_{n,t} - \theta_{m,t})) \quad \forall n \in \mathbf{N}, t \in \mathbf{T} \quad (32)$$

$$Q_{n,t} = V_{n,t} \sum_{m=1}^{N} V_{m,t}((\mathrm{G_{nm}}\sin(\theta_{n,t} - \theta_{m,t}) - (\mathrm{B_{nm}}\cos(\theta_{n,t} - \theta_{m,t})) \quad \forall n \in \mathbf{N}, t \in \mathbf{T} \quad (33)$$

The parameters Conductance $\mathrm{G_{nm}}$ and Susceptance $\mathrm{B_{nm}}$ are the real and imaginary parts of the complex branch Admittance $Y_{nm}$:

$$\mathrm{Y_{nm}} = \mathrm{Z_{nm}^{-1}} = \mathrm{G_{nm}} + \mathrm{jB_{nm}} \quad (34)$$

Note that branch admittance consists of both self-admittances (where $Y_{nm}, n = m$) and mutual admittances (where $Y_{nm}, n \neq m$). The reader is directed to Frank and Rebennack (Frank and Rebennack, 2016) for the computation of complex branch admittances. Note that we assume that shunt susceptance and admittance are negligible, which is a reasonable assumption for short lines or cables, and nominal transformer turns ratio.

In grid-connected DES and microgrids, the main point of connection to the external grid is known as the slack bus, denoted by the subscript $slack$. Active and reactive power injections are not limited at this bus to purport the role of the external grid as both an unlimited source and sink for power, and to obtain feasible solutions for the scenarios tested. The voltage magnitude and angle are fixed at this bus, as shown below:

$$V_{slack,t} = 1, \qquad \theta_{slack,t} = 0 \qquad \forall t \in \mathbf{T} \tag{35}$$

The voltage magnitude at each node and timepoint must remain within pre-specified network bounds:

$$\frac{V^{UB}}{V_{base}} \leq V_{n,t} \leq \frac{V^{LB}}{V_{base}} \qquad \forall n \in \mathbf{N}, t \in \mathbf{T} \tag{36}$$

where $V^{UB}$ and $V^{LB}$ are the per unit values of line voltage upper bound and lower bound, respectively. The voltage angle is also constrained in a similar manner:

$$-180° \leq \theta_{n,t} \leq 180° \qquad \forall n \in \mathbf{N}, t \in \mathbf{T} \tag{37}$$

In a typical OPF formulation, the upper and lower bounds of each generation and storage unit installed are also described. As these are included in the DES formulation, which ultimately gets combined with this OPF formulation, these constraints can be found in Section 2.1.

While the above constraints complete a typical OPF formulation, additional constraints, also known as side constraints, can be included to calculate variables such as the magnitude of current in each branch (Frank and Rebennack, 2016):

$$\left( \left( V_{n,t} \cos \theta_{n,t} - V_{m,t} \cos \theta_{m,t} \right)^2 + \left( V_{n,t} \sin \theta_{n,t} - V_{m,t} \sin \theta_{m,t} \right)^2 \right) * \hat{y}_{nm}^2$$
$$\leq \left( \frac{I^{max}}{S_{base}} \right)^2 \qquad \forall (n,m) \in \mathbf{L}, t \in \mathbf{T} \tag{38}$$

where $I^{max}$ is the line current upper bound defined by the distribution network operator. The bus injection model relies on active and reactive power injections at each node. Nodes that neither consume nor generate power have net active and reactive power injections set to zero, as shown below:

$$P_{n,t} = 0 \quad \forall n \notin \mathbf{M}, t \in \mathbf{T} \tag{39}$$

$$Q_{n,t} = 0 \quad \forall n \notin \mathbf{M}, t \in \mathbf{T} \tag{40}$$

Note that $n \in \mathbf{M} \subseteq \mathbf{N}$ represents the nodes that do generate and consume power, and the constraints for these are described in Section 2.3.

## 2.3.  DES-OPF Linking constraints

The distribution network topology typically indicates which building $i$ is located at which node $n$. Nodes at which buildings are located, that act as either loads or generators, are described as load/generator nodes $n \in \mathbf{M} \subseteq \mathbf{N}$ in this formulation. An indicator $A_{n,i}$ is used to show the connection between such a node $n$ and building $i$:

$$A_{n,i} \in \{0,1\} \quad \forall n \in \mathbf{M}, i \in \mathbf{I} \tag{41}$$

To connect the OPF formulation described in Section 2.2 with the DES formulation summarised in Section 2.1, load/generator nodes, which are potential sites for installing DERs, are assumed to be PQ nodes. The net active and reactive power injections at these nodes are typically expressed using the equations below:

$$P_{n,t} = \sum_{i \in \mathbf{I}} \mathrm{A}_{\mathrm{n,i}} (P_{i,t}^{Gen} - P_{i,t}^{Load}) \quad \forall n \in \mathbf{M}, t \in \mathbf{T} \tag{42}$$

$$Q_{n,t} = \sum_{i \in \mathbf{I}} \mathrm{A}_{\mathrm{n,i}} (Q_{i,t}^{Gen} - Q_{i,t}^{Load}) \quad \forall n \in \mathbf{M}, t \in \mathbf{T} \tag{43}$$

Where $P_{i,t}^{Gen}$ and $Q_{i,t}^{Gen}$ represent the total active and reactive power generated, while $P_{i,t}^{Load}$ and $Q_{i,t}^{Load}$ represent the total reactive power consumed. Eq. (44) is modified to align with the DES baseline and OPF formulations:

$$P_{n,t} = \sum_{i \in \mathbf{I}} \mathrm{A}_{\mathrm{n,i}} \left( \frac{E_{i,t}^{PV,sold} - E_{i,t}^{grid} - \sum_c E_{i,t,c}^{grid,charge}}{S^{base}} \right) \quad \forall n \in \mathbf{M}, t \in \mathbf{T} \tag{44}$$

$$Q_{n,t} = \sum_{i \in \mathbf{I}} \mathrm{A}_{\mathrm{n,i}} \left( \frac{Q_{i,t}^{Gen} - Q_{i,t}^{Load}}{S^{base}} \right) \quad \forall n \in \mathbf{M}, t \in \mathbf{T} \tag{45}$$

where $E_{i,t}^{PV,sold}$ is the excess electricity sold to the grid, $E_{i,t}^{grid}$ is the electricity purchased from the grid to satisfy consumer, and $E_{i,t,c}^{grid,charge}$ is the total electricity purchased to charge the batteries $c \in \mathbf{C}$ installed at each house. Note that $Q_{i,t}^{Load}$ is calculated using a constant power factor applicable to each house, and $Q_{i,t}^{Gen}$ is assumed to be zero in this instance as no reactive power generators are considered.

This completes the linking of the two formulations, and the complete model can be solved using any of the methods described in Section **Error! Reference source not found.**.

## 3. Case Studies

The first case study, labelled Case 1, is a European low voltage (0.4 kV) network supplying electricity to residential customers (Papathanassiou et al., 2005), also used by Morvaj et al. (2016). A diagram of the network including network parameters is presented in Figure 3. It is assumed that all residential consumers are connected as three-phase loads to maintain the assumption that the network is balanced. The power factor (PF) of solar panels is set to unity (i.e., 1.0) as inverters at residential buildings are typically set to provide active power only. The power factor of the loads is 0.85. The case study prescribes an upper bound of 250 A on line current. The peak energy consumption and available space for DER installation at each residential building is recorded in

Table 1.

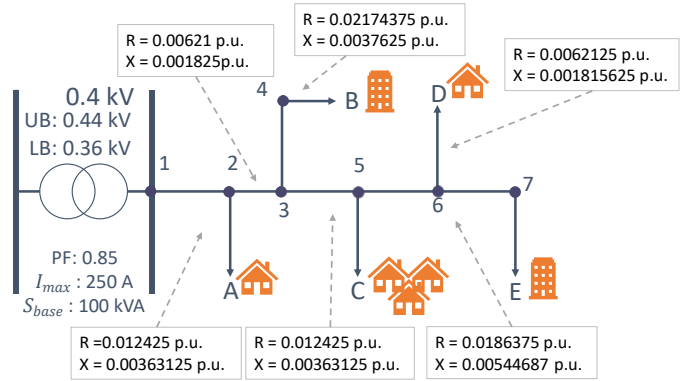

Figure 3. Line diagram of the low voltage distribution network with line parameters Resistance (R) and Reactance (X). The upper and lower bounds for voltage are indicated by UB and LB.

Table 1. Average daily demands for electricity and heating per day

| Building | Peak Electricity (kW) | Peak Heat (kW) | Area available (m²) | Volume available (m³) |
|---|---|---|---|---|
| A | 3.8 | 10.0 | 150 | 5 |
| B | 18.4 | 47.6 | 700 | 5 |
| C | 14.1 | 27.1 | 600 | 5 |
| D | 3.8 | 6.4 | 150 | 5 |
| E | 12.0 | 31.3 | 550 | 5 |
| *Total* | *52.2* | *122.3* | *2150* | *25* |

The second case study, labelled Case 2, is a modified version of the 906-node and 55-load IEEE EU LV network (IEEE, 2020). The modified version contains a transformer with nominal tap ratio (such that phase change is not considered), and includes 22 balanced loads and associated network connections, as shown in Figure 4. The available area for PV installation at each house node is assumed to be 35 m², while the volume available for battery installation is taken as 0.5 m³. Each load has a constant power factor of 0.95, as specified in the case study. It is assumed that all 22 consumers are opting to install renewable energy and storage technologies.

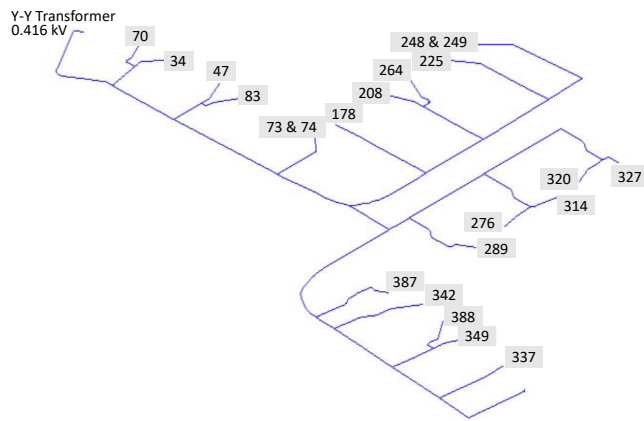

Figure 4. Modified IEEE EU LV network with 22 loads.

The models consider a discretised temporal horizon consisting of 24 hourly timepoints. Averaged daily profiles for electricity and heating demand are used for each season for each case study. Note that the average seasonal electricity profiles for the IEEE case study were derived from the daily 1-min interval profiles provided by the test case, while heating profiles were assumed. Both residential networks are assumed to be located in the UK, therefore averaged solar irradiance profiles (Met Office, 2020) (Figure 5) and the FIT scheme (Ogfem, 2018) applicable to the UK have been used. It is assumed that the installed DES will operate for 20 years in total. This also aligns with the FIT scheme, which promised fixed tariffs for small-scale renewable energy generation and export to the electricity grid, provided that the export occurs only if excess electricity is available after the demand of the consumer has been satisfied. In addition to these tariffs, the Economy 7 pricing strategy is considered (The Consumer Council, 2018), which provides different day and night electricity prices to consumers in the UK. As the day price is typically higher than the night price, this assists the study of battery installation and operation within the DES.

Other key input parameters associated with the DERs and the system are all provided in Appendix A. Python/Pyomo-based models (Hart et al., 2011), all input data used, and results files have been made available via Github: https://github.com/Ishanki/DES-OPF-Design.

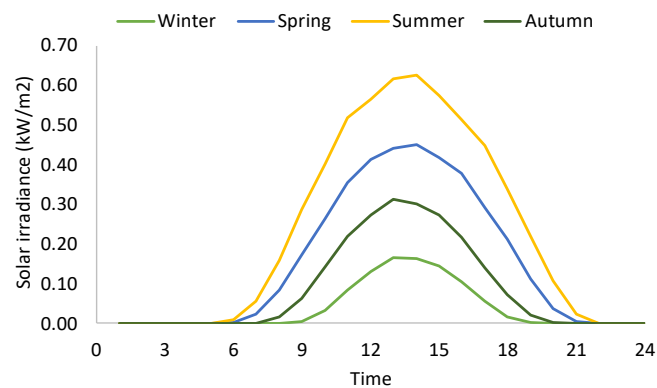

Figure 5. Averaged daily solar irradiance profiles for each season.

## 4.  Results and Discussion

The MILP master problems have been solved using CPLEX (IBM, 2019), while the bi-level nonlinear subproblems have been solved using CONOPT (Drud, 1985). The MINLP model has been solved using SBB (Bussieck and Drud, 2001). Both case studies were tested using an Intel® Core™ i7-10510U CPU at 1.80G - 2.30 GHz. Case 1 MILP contains 5,693 continuous variables, 980 binary variables, and 10,734 constraints. Case 1 DES-OPF models contain a total of 8669 continuous variables, 14,842 constraints, and the same number of binary variables as the MILP (which are either fixed or free, depending on whether BL or MINLP is solved). Case 2 MILP consists of 68,269 continuous variables, 4,312 binary variables, and 134,135 equality and inequality constraints. Case 2 DES-OPF models consists of 168,109 continuous variables, 258,271 constraints, and the same number of binary variables as its MILP. Note that an initialisation strategy was not used in either case study when solving the NLP, to test the capabilities of the commercial solver to find feasible solutions.

### 4.1.  Results

Case 1 results for MILP, BL, and MINLP are summarised in Table 2. As expected, the MILP model produces the best objective value compared to the other models. While no significant differences in objective function are observed and the same PV capacity has been installed across the three models, differences in battery capacity are observed. Further investigation reveals that the higher battery capacity in BL and MINLP models offset the current violations that would otherwise exist at peak PV energy production if the MILP design was implemented, which are portrayed in Figure 6. This is also reflected by the reduced export income, as BL and MINLP both opt for greater electricity storage as opposed to export. Testing the MILP design using modified BL (MILP Check) demonstrates that the MILP design can be feasible, provided that the operational schedule is modified to evade the current violations. This results in a minor loss of income from PV generation and export due to renewable energy curtailment, which, in turn, is reflected in the slightly higher objective value seen in MILP Check. No voltage violations are observed in Case 1.

Table 2. Results for Case 1 with a breakdown of annualised costs.

| Breakdown | MILP | MILP Check | BL | MINLP |
|---|---|---|---|---|
| Time taken (s) | 3.20 | 5.62 | 5.93 | 10.10 |
| Objective value (£) | 32,615 | 33,106 | 33,010 | 32,903 |
| % Difference with MILP obj. | - | 2 | 1 | 1 |
| Relative optimality gap | 0 | - | - | 0 |
| PV investment (£) | 54,231 | 54,231 | 54,231 | 54,231 |
| Boiler investment (£) | 1,440 | 1,440 | 1,440 | 1,440 |
| Battery investment (£) | 585 | 585 | 1,502 | 2,341 |
| Grid electricity (£) | 23,590 | 23,590 | 21,943 | 20,339 |
| PV operation (£) | 3,839 | 3,839 | 3,839 | 3,839 |
| Boiler operation (£) | 5,240 | 5,240 | 5,240 | 5,240 |
| Generation income (£) | 42,855 | 42,529 | 42,674 | 42,791 |
| Export income (£) | 13,697 | 13,533 | 13,134 | 12,708 |

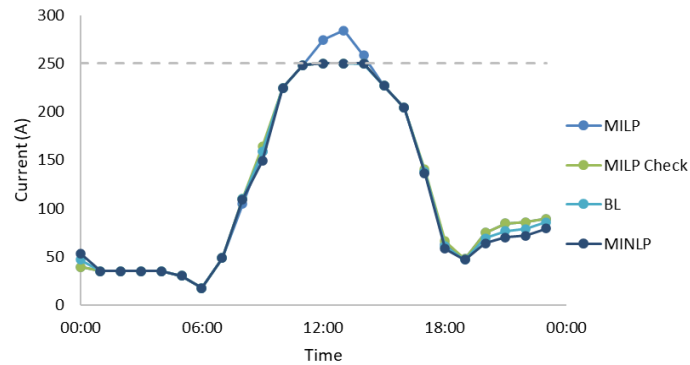

Figure 6. Currents in branch (1, 2) during summer, as calculated in all four models for Case 1. The dashed line represents the line current upper bound of 250 A.

Case 2 results are summarised in Table 3, where much more prominent differences in objectives and designs between the models are observed when compared to Case 1. Note that MINLP results are not included as the solver was unable to find a feasible binary topology prior to exceeding solver time limit set to 3600 s. Although the MILP continues to produce the best objective value, MILP Check reveals that the true objective could be 37% higher than initially predicted, when more detailed network constraints are considered. This is primarily due to the existence of voltage violations causing renewable energy curtailment, which the DC approximations are unable to detect due to the assumption that the network maintains nominal voltage at all nodes. This is portrayed in Figure 7, where voltage violations are observed throughout much of the PV energy production period during summer. Note that an upper bound for line current has not been imposed in this case study. BL attempts to mitigate some of these negative impacts by reducing PV capacity and nearly quadrupling the total battery capacity, utilising the stored energy to reduce electricity purchasing costs. Thus, unlike in Case 1, Case 2 demonstrates that implementing the MILP design in practice may significantly increase renewable energy wastage and overall costs, making this design unsuitable in comparison to that proposed by BL. Case 2 also demonstrates the scalability of BL to find feasible solutions, despite the increase in computational time, when compared to the MINLP which has been unable to find a solution.

Table 3. Results for Case 2 with a breakdown of annualised costs.

| Breakdown | MILP | MILP Check | BL |
|---|---|---|---|
| Time taken (s) | 81.24 | 375.32 | 402.23 |
| Objective value (£) | 21,297 | 29,250 | 27,882 |
| % Difference with MILP obj. | - | 37 | 31 |
| Relative optimality gap | 0 | - | - |
| PV investment (£) | 19,422 | 19,422 | 17,364 |
| Boiler investment (£) | 3,359 | 3,359 | 3,359 |
| Battery investment (£) | 64 | 64 | 253 |
| Grid electricity (£) | 4,962 | 4,962 | 4,764 |
| PV operation (£) | 1,375 | 1,375 | 1,229 |
| Boiler operation (£) | 13,804 | 13,804 | 13,804 |
| Generation income (£) | 15,348 | 10,057 | 9,584 |
| Export income (£) | 6,368 | 3,707 | 3,413 |

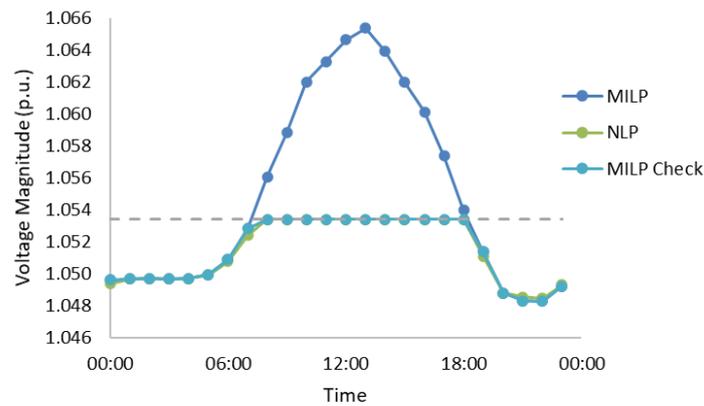

Figure 7. Phase voltage of Bus 388 in summer, as predicted by three out of the four models. The dashed line represents the voltage upper bound. The true voltages for the MILP observed here have been obtained using the AC OPF formulation, as the DC approximations assume nominal voltage.

### 4.2. Verification of power flow formulations

To compare the accuracy of the nonlinear balanced AC power flow formulations presented in Section **Error! Reference source not found.**, the Newton-Raphson algorithm with timeseries power flow analysis on Pandapower (Thurner et al., 2018) is used. BL is used to simulate power flow analysis with the OPF formulation, without the capabilities of installing DES (i.e., power must be purchased from the external grid). Note that Pandapower reports current values as phase currents as opposed to line currents (Thurner et al., 2021), therefore these values have been multiplied by $\sqrt{3}$ to obtain line currents for comparison. Percentage errors for voltage magnitude ($V_{ph}$) and angle ($\theta_{ph}$) at each bus, active ($P_{bus}$) and reactive ($Q_{bus}$) power at each bus, and branch current ($I_L$) are reported in Table 4. As all errors are within ±0.23%, the power flow formulation utilised is deemed sufficiently accurate.

Table 4. Percentage errors calculated between Pandapower solutions and power flow formulations used in this study.

| | $V_{ph}$ | $\theta_{ph}$ | $P_{bus}$ | $Q_{bus}$ | $I_L$ |
|---|---|---|---|---|---|
| Min | 0.000 | -0.029 | -0.211 | -0.036 | -0.065 |
| Max | 0.000 | 0.018 | 0.221 | 0.029 | 0.063 |

## 4.3.    Discussion and Limitations

Returning to the assessment of practical feasibility of the resulting designs, Table 5 compares the results of the models for both case studies. The operational schedule proposed by the MILP (representative of DES baseline models) for both cases, would be infeasible with respect to the underlying network. However, tests conducted with MILP Check confirmed that these designs would be feasible, provided that excess renewable energy is curtailed during normal operation due to inadequate battery storage. In Case 1, these losses did not have a significant impact on the overall annualised costs for the consumer. However, in Case 2, a significant impact on annualised costs was observed due to these losses. Despite the variability of the impacts on costs, these results confirm that DES baseline models are ill-suited for designing grid-connected DES. It is evident that a DES-OPF framework is essential to accurately represent the underlying network and its constraints and adjust the DES design to minimise energy losses and maximise energy utilisation for baseline operation. While these results have been observed for baseline operation, it can be conjectured that energy losses and network violations observed here could be amplified in the presence of disturbances or unexpected changes to operation, resulting in worse performance and financial losses. The results also emphasise the important role of energy storage, acting as a strategic buffer between DES operation and meeting network requirements, while maximising local renewable energy consumption. This suggests that small-scale batteries can play a significant role in minimising the negative impacts to distribution network infrastructure due to increasing small-scale renewables. This is an interesting finding as small-scale PV deployment (0 − 4 kW) exceeded 600,000 installations (Nolden, 2015) in 2015, providing over 3000 GW (BEIS, 2021), while behind-the-meter battery capacity has been slowly rising to a mere 7 GW in 2020 (IEA, 2021). With battery investment often taking a backseat in comparison to renewable generation investment, especially in terms of available grants and tariffs for small-scale DES prosumers in the U.K., policy makers interested in reducing carbon emissions and increasing local renewables could endeavour to make energy storage more accessible and affordable.

Table 5. Assessing the practical feasibility of the resulting designs for Case 1 and 2. The asterisk (*) indicates the use of MILP Check to confirm whether network constraints are violated.

| Case Study | Model | Network constraints respected? | Reduces energy wastage? |
|---|---|---|---|
| Case 1 | MILP | Yes* | No |
| | BL | Yes | Yes |
| | MINLP | Yes | Yes |
| | | | |
| Case 2 | MILP | Yes* | No |
| | BL | Yes | Yes |

The bi-level method proposed in this study to solve the DES-OPF design framework is not without limitations. Firstly, there may be instances where the fixed binary topology needs to be updated in order to avoid infeasibility or find better feasible solutions. In both test cases, BL has returned feasible solutions as it contains sufficient degrees of freedom provided by continuous design and operational variables. Note that the MINLP performs better than BL in Case 1 as it has the greatest degrees of freedom, due to the retention of free binary variables along with all continuous variables. The BL method can be further improved by incorporating a feedback system to Level 1 (MILP), should Level 2 (NLP) render infeasible. Secondly, we utilise balanced OPF formulations, as unbalanced multiphase formulations are much more complex to formulate and solve. As evident in Case 2 results, the balanced formulation successfully detects voltage violations observed in the original unbalanced network. Thus, we have deemed that the balanced approximation provides the accuracy-complexity balance needed to solve larger test cases. However, in future work, we aim to develop algorithms capable of solving the more complex DES-Unbalanced OPF problem and understand how incorporating an unbalanced formulation may impact the design and objective.

Overall, the results demonstrate that the bi-level method proposed in this study is capable of finding feasible solutions to larger problems using mathematical programming approaches, while maintaining an accuracy-complexity balance, as opposed to previous approaches in literature. The results also shed light on the synergies between DES and the underlying distribution network, which have been previously overlooked in research, and highlight the importance of considering network configuration within DES frameworks.

## 5. Conclusions

This study aims to assess the practical feasibility of DES designs proposed by conventional MILP frameworks that utilise oversimplified DC approximations, in comparison with higher-fidelity DES-OPF frameworks. A new bi-level method for solving the nonconvex DES-OPF problem is proposed, which includes the MILP framework with DC approximations at Level 1, and a complete nonconvex OPF formulation at Level 2. The method proposed avoids the use of iterative numerical checks that require the use of external tools, unscalable linearisations, and unreliable and computationally expensive metaheuristic approaches for solving the overarching MINLP. The method prioritises feasibility and accuracy of the solution over global optimality. Two test cases of varying scale are utilised to test the proposed method, and to compare designs between the DES MILP and combined DES-OPF frameworks. Resulting designs from the MILP model are tested with a modified version of the bi-level method (containing detailed OPF) to assess their practical feasibility. Both test cases show that the MILP with DC approximations initially predicts the lowest total annualised costs in the absence of detailed power flow constraints, but produce the highest costs when the designs are tested using the modified bi-level method. This is due to energy curtailment, where its impacts are observed most prominently in Case 2, resulting in much higher costs than initially predicted. On the other hand, the bi-level method produces practically feasible models by employing greater battery storage to minimise renewable energy curtailment and network violations during baseline operation. Case 2, which has a greater number of buses and loads compared to Case 1, demonstrates that the proposed method is capable of finding feasible solutions when compared with more computationally expensive MINLP methods that have been employed in previous studies. Overall, the study affirms that the combined DES-OPF frameworks are essential when designing grid-connected DES. Future work involves further development of the bi-level method to solve unbalanced multiphase power flow while maintaining an accuracy-complexity balance. Such efforts can further confirm the effectiveness of using more detailed optimisation models to design and operate robust DES while protecting existing distribution network infrastructure.

## Acknowledgement

Ishanki De Mel would like to thank Dr. Boran Morvaj and Dr. Ralph Evins for kindly sharing the demand data used in their article, ref: Morvaj et al. (2016).

## Appendix A – Input Parameters

Table A 1. Technology- and pricing-related parameters used in the models

| Parameter | Value | Units | Reference |
|---|---|---|---|
| Interest rate | 0.075 | | Assumed |
| DES lifetime | 20 | | (Mariaud et al., 2017) |
| Economy7 night tariff | 0.08 | £/kWh | (Department for Business, 2021) |
| Economy7 day tariff | 0.18 | £/kWh | (Department for Business, 2021) |
| Gas purchasing price | 0.02514 | £/kWh | Private correspondence[1] |
| FIT export tariff | 0.0503 | £/kWh | (Ofgem, 2017) |
| FIT generation tariff | 0.1 | £/kWh | Assumed[2] |
| **PVs** | | | |
| PV Investment cost | 450 | £/panel | Average from (Clissit, 2020) |
| PV Efficiency | 0.18 | | (Mariaud et al., 2017) |
| PV Fixed operational cost | 12.5 | £/kW-yr | (Mariaud et al., 2017) |
| Panel area | 1.75 | $m^2$ | Approximation from (Clissit, 2020) |
| Panel capacity | 0.25 | kW | Average from (Clissit, 2020) |
| **Boilers** | | | |
| Boiler investment cost | 120 | £/kW | (Myers et al., 2018) |
| Boiler efficiency | 0.94 | | (Alsop, 2020) |
| **Batteries** | | | |
| Volumetric energy density | 20 | $kWh/m^3$ | (Mariaud et al., 2017) |
| Max DoD | 0.85 | | (Mariaud et al., 2017) |
| Max SoC | 0.9 | | (Mariaud et al., 2017) |
| Investment cost | 270 | £/kWh | (Mariaud et al., 2017) |
| Operational cost | 11 | £/kWh-yr | (Mariaud et al., 2017) |
| Round trip efficiency (RTE) | 0.89 | - | (Mariaud et al., 2017) |
| Charge efficiency | 0.94 | - | Calculated using RTE |
| Discharge efficiency | 0.91 | - | Calculated using RTE |

---

[1] The University of Surrey Estates and Facilities Department staff
[2] The FIT generation tariff from ref. (Ofgem, 2017) was increased to encourage PV installation.